\documentclass[twoside,11pt,a4 paper]{article}
\usepackage{amssymb,amsmath,latexsym,theorem}
\usepackage{graphics}
\usepackage{amssymb}
\usepackage{mathrsfs}
\usepackage{makeidx}
\usepackage{color}
\usepackage{eucal}
\theoremheaderfont{\bf}
\newtheorem{pro}{Proposition}[section]
\newtheorem{teo}[pro]{Theorem}

\newtheorem{ob}[pro]{Remark}
 \newtheorem{cor}[pro]{Corollary}

\newtheorem{defin}[pro]{Definition}

\newcommand{\pp}{\partial}

\title{\bf Further properties of the Bergman spaces of slice regular functions}

\author{Fabrizio Colombo\\Politecnico di
Milano\\Dipartimento di Matematica\\Via E. Bonardi, 9\\20133 Milano,
Italy\\fabrizio.colombo@polimi.it
\and
J. Oscar Gonz\'alez-Cervantes\\
 Departamento de Matem\'aticas  \\
 E.S.F.M. del
I.P.N. 07338 \\
M\'exico D.F., M\'exico
\\
jogc200678@gmail.com
\and
Irene Sabadini\\Politecnico di
Milano\\Dipartimento di Matematica\\Via E. Bonardi, 9\\20133 Milano,
Italy\\irene.sabadini@polimi.it }
\date{  }
\begin{document}
\maketitle

\begin{abstract}
In this paper we continue the study of Bergman theory for the class of slice regular functions.
  In  the  slice regular setting there are two possibilities to introduce the Bergman spaces, that are called  of  the  first and of the second kind.
In this paper we mainly consider the Bergman theory of the second kind, by providing an explicit description of the Bergman kernel in the case of the unit ball and of the half space. In the case of the unit ball, we study the Bergman-Sce transform. We also show that the two Bergman theories can be compared only if suitable weights are taken into account. Finally, we use the Schwarz reflection principle to relate the Bergman kernel with its values on a complex half plane.
\end{abstract}



\section{Introduction }
The literature on Bergman theory is wide and
as classical  reference books we mention, with no claim of completeness, the books \cite{berg}, \cite{bgsc} and \cite{K}.
The theory has also been developed for hyperholomorphic functions such as the quaternionic regular  functions in the sense of Fueter
and the theory of monogenic functions,  \cite{bradel}, \cite{bds},
  \cite{const},  \cite{conskrau}, \cite{conskrau2}, \cite{delanghe}, \cite{shavas}, \cite{shavas2}, \cite{shavas3}  and the literature therein.
  Recently we have started the study of Bergman theory in the slice hyperholomorphic setting,
  see \cite{CGLSS}, \cite{CGS1}, \cite{CGS2}, \cite{CGS3}.
  In this framework, it is possible to give two different notions of Bergman spaces, the so-called Bergman spaces of the first and of the second kind. Bergman spaces of the first kind in the slice regular setting are mainly studied in \cite{CGLSS}. They are defined as
$$ \mathcal A(\Omega):=\{  f\in \mathcal{SR}(\Omega) \ |\    \displaystyle \|f\|^2_{\mathcal A(\Omega)}:= \int_{\Omega}|f|^2 d\mu <\infty \}$$
where $\Omega$ is an axially symmetric bounded open set in the space of quaternions $\mathbb H$ and $\mathcal{SR}(\Omega)$ denotes the set of slice regular functions on $\Omega$.
 The Riesz representation theorem allows to write the so--called slice regular  Bergman  kernel  of  the  first  kind  $\mathcal B(\cdot,\cdot)$ associated with $\Omega$, leading  to  the  integral  representation
 \begin{equation}\label{Bergmanreproduction}
f(q)=\int_{\Omega} \mathcal B(q,\cdot) f d\mu,\quad \forall f\in \mathcal A(\Omega).
\end{equation}

If we denote by $ \mathbb S^2$ the sphere of pure imaginary quaternions and by $\mathbb{C}({\bf i})$ the complex plane with imaginary unit
${\bf i}\in \mathbb S^2$, then slice regular functions $f: \Omega \to\mathbb{H}$ are functions whose restrictions to  $\Omega_{{\bf i}}:=\Omega\cap\mathbb C({\bf i})$, for
every $ {\bf i}\in \mathbb S^2$, are holomorphic maps. For a study of this class of functions see the books \cite{bookfunctional}, \cite{GenSS} and the references therein.
 When we suppose that the domain $\Omega $ intersects the real line and it is axially symmetric, i.e. symmetric with respect to the real axis, then the restriction of a slice regular function to two different planes, $ \mathbb C({\bf i})$ and  $\mathbb C({\bf j})$ with ${\bf i}\not={\bf j}$, turn out to be strictly related by the Representation Formula
 $$
 f(x+y{\bf{i}})=\frac{1}{2}(1-{\bf i} {\bf j}) f(x+y{\bf j})+ \frac{1}{2}(1+{\bf i} {\bf j}) f(x-y{\bf j})
 $$
which asserts that if we know the value of the function $f$ on the complex plane $\mathbb C({\bf j})$ then we can reconstruct $f$ in all the other planes $\mathbb C({\bf i})$ for every  ${\bf{i}}\in \mathbb{S}^2$.
So the global behavior of these functions on these particular axially symmetric sets
is in fact completely  determined by their behavior on a complex plane $\mathbb C({\bf i})$.

From this fact, it follows that there is a second way to define the Bergman space and the Bergman kernel: we can work with the restriction to a complex plane $\mathbb C({\bf i})$ and then to extend it using the Representation Formula. This approach gives rise  the so-called Bergman  theory  of  the  second  kind.

In this paper we will  mainly work with Bergman spaces of the second kind, see Section 3. In section 4, we will provide the explicit description of the Bergman kernel in the case of the unit ball and of the half space. The two Bergman theories of the first and of the second kind are compared in Section 5, where we show that they contain different elements, unless some suitable weights are taken into account.
Then, in Section 6, we study the
 Bergman-Fueter transform which is an integral transform that associates to every slice regular function $f$ defined on $\Omega$ the Fueter regular function $\breve{f}$ given by $\breve{f}=\Delta f$, where $\Delta$ is the Laplace operator.
Finally, in Section 7, we use the Schwarz reflection principle to relate the Bergman kernel of the second kind with its values
on the half plane $\mathbb C^+({\bf i})=\{x+y{\bf i}, \ y\geq 0\}$ and also to write the inner product of two elements using an integral on an half plane.

\medskip
\section{Preliminary results}
By $\mathbb H$ we denote the algebra of real quaternions.
A quaternion $q$ is an element of the form $q=q_0+{\imath} q_1+{\jmath} q_2+{k} q_3$ where $i,j,k$ satisfy $\imath^2=\jmath^2=k^2=-1$ and $\imath\jmath=-\jmath\imath=k$, $\jmath k=-k\jmath=\imath$, $\imath k=-k\imath=\jmath$. A quaternion will also be denoted as
 $q=q_0+\underline{q}$ where $q_0$ and $\underline{q}=\imath q_1+\jmath q_2+kq_3$ are its real and imaginary part, respectively. The modulus of a quaternion $q$ is defined as $|q|=(q_0^2+q_1^2+q_2^2+q_3^2)^{1/2}$.
 Let
$$\mathbb S^2:=\{ {{\underline{q}}=\imath q_1+\jmath q_2+kq_3}  \ \mid  \  | {\underline{q}} |_{\mathbb R^3}=1 \}.$$
It is important to note that an element ${\bf i}\in\mathbb{S}^2$ is such that  ${\bf i}^2=-1$.
Given   ${\bf i}\in\mathbb{S}^2$,  we denote by $\mathbb C({\bf i})$   the real linear space generated by $1$ and ${\bf i}$.
It is immediate that $\mathbb C({\bf i})  \cong \mathbb C$.  \\
 For any open set $\Omega\subset\mathbb H$ let
 $$
\Omega_{{\bf i}}:=\Omega\cap\mathbb C({\bf i}), \quad
\Omega^+_{{\bf i}}= \{x+y{\bf i}\in \Omega_{{\bf i}}\ | \ y\geq 0\}.
$$
We set $\mathbb B^4:=\{ {q} \in \mathbb H \ \mid   \  |  q | <1 \}$, and
$\mathbb{D}_{{\bf i}}:=\mathbb{B}^4\cap\mathbb{C}(\mathbf{i})$.
 Finally, for any function  $f: \Omega \to \mathbb{H}$   we denote its restriction to $\Omega\cap\mathbb C({\bf i})$  by $f_{\mid_{\Omega_{{\bf i}}}}$.
\\
Let us now recall the definition of slice regular functions.
\begin{defin}
A real differentiable quaternion-valued function $f$ defined on an open set $\Omega\subset\mathbb{H}$ is called (left) slice regular on $\Omega$ if  for any ${\bf i}\in \mathbb S^2$ the function $f_{\mid_{\Omega_{{\bf i}}}}$
 is such that
$$
 \left( \frac{\partial}{\partial x} + {\bf i}  \frac{\partial}{\partial y}\right) f_{\mid_{\Omega_{{\bf i}}}}(x+y{\bf i})=0 \textrm{ on $\Omega_{{\bf i}}$}.
 $$
 We denote by $\mathcal{SR}(\Omega)$ the set of (left) slice regular functions on $\Omega$.\\
The function $f$ is called slice anti-regular on the right if for any ${\bf i}\in \mathbb S^2$ the function $f_{\mid_{\Omega_{{\bf i}}}}$  is such that
$$
f_{\mid_{\Omega_{{\bf i}}}}(x+y{\bf i})\left( \frac{\partial  }{\partial x}- {\bf{i}}\frac{\partial}{\partial y}\right)=0 \textrm{ on $\Omega_{{\bf i}}$}.
$$
\end{defin}
\begin{ob} {\rm  With suitable modifications, we can define slice right regular functions and slice left anti-regular functions.
}
\end{ob}
We now introduce a special subclass of domains on which slice regular functions have nice properties.
\begin{defin}
Let $\Omega \subseteq \mathbb{H}$ be a domain. We say that $\Omega$ is a
\textnormal{slice domain (s-domain} for short) if  {} $\Omega \cap \mathbb{R}$ is non empty and
if {} $ \Omega\cap\mathbb C({\bf i})$ is a domain in $\mathbb C({\bf i})$ for all ${\bf i} \in \mathbb{S}^2$.
\end{defin}

\begin{defin}\label{def_circular}
Let $\Omega \subseteq \mathbb{H}$. We say that $\Omega$ is
\textnormal{axially symmetric} if  for every $x+{\bf i}y \in \Omega$ we have that
$x+{\bf j}y\in\Omega$ for all ${\bf{j}}\in \mathbb{S}^2$.
\end{defin}
Note that  any axially symmetric open set $\Omega$ can be uniquely associated with an open set $O_\Omega\subseteq\mathbb{R}^2$, the so-called {\em generating set}, defined by
$$
O_\Omega:=\{(x,y)\in\mathbb{R}^2\ |\ x+y\mathbf{i}\in\Omega,\ \mathbf{i}\in\mathbb{S}^2\}.
$$
In the sequel, we will denote by $Hol(\Omega_{\bf i})$ the set of $\mathbb C_{\bf i}$-valued functions, holomorphic on $\Omega_{\bf i}\subseteq \mathbb C_{\bf i}$.
Let $Z$ denote the complex conjugation, i.e.  $Z(z)=\bar z,\quad \forall z\in \mathbb C({\bf i})$.
If $\Omega$ is an axially symmetric s-domain then for any ${\bf i} \in \mathbb S^2$ the domain $\Omega_{\bf i}\subset \mathbb C({\bf i})$  satisfies that $Z(\Omega_{\bf i})=\Omega_{\bf i}$.
We also define $$Hol_c(\Omega_{\bf i}):=\{f\in Hol(\Omega_{\bf i}) \ \mid \ Z\circ f\circ Z=f\}.$$
Functions belonging to $Hol_c(\Omega)$ are called in the literature {\em real} or {\em intrinsic}.

Slice regular functions satisfy the following property (see \cite{CGS, newadvances, csTrends}).
\begin{teo}[Representation Formula]\label{formula} \index{Representation Formula} Let
$f$ be a slice regular function on an axially symmetric s-domain $\Omega\subseteq  \mathbb{H}$. Choose any
${\bf{j}}\in \mathbb{S}^2$.  Then the following equality holds for all $q=x+{\bf{i}}y \in \Omega$:
\begin{equation}\label{relationflats}
    f(x+y{\bf{i}})=\frac{1}{2}(1-{\bf i} {\bf j}) f(x+y{\bf j})+ \frac{1}{2}(1+{\bf i} {\bf j}) f(x-y{\bf j}).
     \end{equation}
\end{teo}
\begin{ob}{\rm By the Splitting Lemma, the restriction of a slice regular function to any complex plane $\mathbb{C}({\bf i})$ can be written as $f_{|\Omega_{\bf i}}(z)=F(z)+G(z) {\bf j}$ where ${\bf j}$ is an element in $\mathbb S^2$ orthogonal to ${\bf i}$ and $F,G:\Omega_{\bf i}\to\mathbb C({\bf i})$ are two holomorphic functions. Conversely, the Representation formula allows to extend a function of the form $f(z)=F(z)+G(z) {\bf j}$, defined on $\Omega_{\bf i}$ to a slice regular function defined on $\Omega$. Thus we have the so-called extension operator
$$ P_{\bf i} :   Hol(\Omega_{\bf i})+ Hol(\Omega_{\bf i}){\bf j} \longrightarrow  \mathcal{SR}(\Omega),$$ defined for any $f\in  Hol(\Omega_{\bf i})+ Hol(\Omega_{\bf i}){\bf j}$  by
 \begin{equation}\label{extPi}
   P_{\bf i}[f](q)=P_{\bf i}[f](x+I_qy)=\frac{1}{2}\left[(1+ I_q{\bf i})f(x-{\bf i}y) + (1- I_q {\bf i}) f(x+ {\bf i}y)\right],
  \end{equation}
  where $q=x+I_qy$.
  }
  \end{ob}

\subsection{The slice regular Bergman theory of the first kind}

Let  $\Omega\subset\mathbb H$ be a bounded axially symmetric s-domain. We will consider the space  $\mathcal L_2(\Omega,\mathbb H)$ formed by  functions $f:\Omega\to  \mathbb H$ such that
$$
\int_{\Omega}|f|^2 d\mu < \infty,
$$
where $d \mu$ denotes the Lebesgue volume element in $\mathbb R^4$.
The functional $\langle\cdot, \cdot \rangle_{
\mathcal L_2(\Omega,\mathbb H)}:\mathcal L_2(\Omega,\mathbb H)\times \mathcal L_2(\Omega,\mathbb H) \to \mathbb H$, defined by
 $$
 \langle f, g \rangle_{\mathcal L_2(\Omega,\mathbb H)}=\int_{\Omega} \bar f g d\mu,
 $$
is an $\mathbb H$-valued inner product on $\mathcal L_2(\Omega,\mathbb H)$.
The space $\mathcal L_2(\Omega,\mathbb H)$  is then equipped with the norm
\begin{equation}\label{normaglob}
\|f \|_{\mathcal L_2(\Omega,\mathbb H)}:= \left(\int_{\Omega}|f|^2 d\mu\right)^{\frac{1}{2}}.
\end{equation}

%
%

\begin{defin}\label{Bergmanfirtsk}  The set   $\mathcal A(\Omega
):=  \mathcal{SR}(\Omega)\cap \mathcal L_2(\Omega,\mathbb H)$, equipped with the norm and the inner product inherited from  $\mathcal L_2(\Omega,\mathbb H)$, is called the slice regular Bergman space of the first kind associated with $\Omega$.
\end{defin}
The following results have been proved in \cite{CGLSS}, see Proposition 3.4 and Theorem 3.5, but we repeat them since they will be useful in the sequel.

\begin{pro}\label{compactsetR4}Let $\Omega$ be a bounded axially symmetric s-domain.
For any compact set $K\subset  \ \Omega \setminus \mathbb R$,  there exists $\lambda_K>0$ such that
$$\sup\{ |f(q)| \ \mid \ q\in K  \}\leq \lambda_K\| f\|_{\mathcal A(\Omega)},\quad \forall f\in  \mathcal A(\Omega) .  $$
\end{pro}
\begin{teo}\label{completeness} The set $\mathcal A(\Omega)$ is a complete  space. \end{teo}

 Theorem \ref{completeness} implies that
$ \mathcal A(\Omega)$ is a quaternionic right linear Hilbert space.

For further properties on the  Bergman theory of the first kind we refer the reader to \cite{CGLSS}.

\section{The slice regular Bergman theory of the second kind}

 Let $\Omega\subseteq  \mathbb{H}$ be a bounded axially symmetric s-domain.
We introduce here a family of Bergman-type spaces. For every ${\bf{i}}\in \mathbb{S}^2$  we  set
  $$
\displaystyle{     \mathcal A  (\Omega_{ \bf i  } ):=    \left\{    f\in \mathcal{SR}(\Omega) \   \mid  \
    \displaystyle \|f\|^2_{{\mathcal A}(\Omega_{\bf i})}:= \int_{\Omega_{\bf i}}|f_{\mid_{\Omega_{{\bf i}}}}|^2 d\sigma_{\bf i} <\infty   \right\} },
  $$
  where $d\sigma_{{\bf i}}$ denotes the Lebesgue measure in the plane $\mathbb C_{\bf i}$ and $d\sigma_{\bf i}(x+{\bf i}y)=dx dy$.
On $\mathcal A(\Omega_{{\bf i}})$ we define the norm
\begin{equation}\label{normSecondBergman}
  \displaystyle\|f\|_{{\mathcal A}(\Omega_{\bf i})}:= \left(\int_{\Omega_{\bf i}}|f_{\mid_{\Omega_{{\bf i}}}}|^2 d\sigma_{\bf i}\right)^{\frac{1}{2}}, \quad \forall f\in \mathcal A(\Omega_{\bf i}),
\end{equation}
and  the inner product
\begin{equation}\label{scalar}
\langle f,g\rangle_{\mathcal A(\Omega_{{\bf i}})}=\int_{\Omega_{{\bf i}}} \overline{f}\, g\,d\sigma_{{\bf i}}.
\end{equation}

\begin{pro}\label{pro1} Let ${\bf i}, {\bf j}$ be any unit  vectors, and let $\Omega$ be a bounded axially symmetric s-domain and $f\in\mathcal {SR}(\Omega)$. Then $ f\in  \mathcal A(\Omega_{\bf i})$
 if and only if  $ f\in  \mathcal A(\Omega_{\bf j})$.
   \end{pro}
Proof.
 Using formula (\ref{relationflats}), we have that
  $$
  |f(x+y{\bf{i}})|^2\leq | f(x+y{\bf j})|^2 + 2 |f(x+y{\bf j})||f(x-y{\bf j})|+ |f(x-y{\bf j})|^2.
  $$
As
  $$
  |f(x+y{\bf j})||f(x-y{\bf j})|\leq \frac{1}{2} (|f(x+y{\bf j})|^2+ |f(x-y{\bf j})|^2),
  $$
then
  \begin{equation}
    \label{inequalitymodules}
    |f(x+y{\bf{i}})|^2\leq 2 \left[ | f(x+y{\bf j})|^2 + | f(x-y{\bf j})|^2\right],
  \end{equation}
 and  integrating both sides of (\ref{inequalitymodules}) on the generating set $O_\Omega$, one obtains that
$$\int_{{\Omega_{\bf i}}} |f_{|{\Omega_{\bf i}}}   |^2d\sigma_{\bf i} \leq 2  \int_{{\Omega_{\bf j}}} |f_{| {\Omega_{\bf j}}}   |^2d\sigma_{\bf j}  . $$
Now replacing ${\bf i}$ by ${\bf j}$ in \eqref{inequalitymodules}, we get
$$\int_{{\Omega_{\bf j}}} |f_{| {\Omega_{\bf j}}}   |^2d\sigma_{\bf j} \leq 2  \int_{{\Omega_{\bf i}}} |f_{| {\Omega_{\bf i}}}   |^2d\sigma_{\bf i}, $$
and the assertion follows.
\hfill $\blacksquare$
\begin{cor} Given any   ${\bf i}, {\bf j} \in \mathbb S^2 $, the Bergman spaces   $ \mathcal A(\Omega_{\bf i})$  and  $ \mathcal A(\Omega_{\bf j}) $  contain the same elements and have equivalent norms.
\end{cor}

We also recall the following result, see \cite{CGLSS}:
\begin{teo}[Completeness of $\mathcal A(\Omega_{\bf i})$]\label{complfOi}
 Let $\Omega$ be a bounded axially symmetric s-domain. The spaces $\left( \mathcal A(\Omega_{\bf i}), \|\cdot\|_{{\mathcal A}(\Omega_{\bf i})} \right)$ are complete for every ${\bf i} \in \mathbb S^2$.
  \end{teo}

\begin{ob}\label{sliceBergmankernel}
{\rm
 Let  $\Omega$ be a given bounded axially symmetric s-domain in $\mathbb{H}$.
 Theorem \ref{complfOi} implies that the  Bergman spaces $\mathcal A(\Omega_{\bf i})$, for ${\bf i} \in \mathbb S^2$, are  quaternionic right linear Hilbert spaces.
}
\end{ob}
\begin{defin}
We say that a function $f\in \mathcal S\mathcal R(\Omega)$ belongs to the Bergman space of the second kind $\mathscr A(\Omega)$ on the bounded, axially symmetric s-domain $\Omega$ if $f\in\mathcal A(\Omega_{\bf i})$ for some ${\bf i}\in\mathbb S^2$.
\end{defin}
The norm $\|\cdot\|_{\mathscr A(\Omega)}$ in $\mathscr A(\Omega)$ is defined to be one of the equivalent norms $\|\cdot\|_{{\mathcal A}(\Omega_{\bf i})}$ for some fixed $\bf i\in\mathbb S^2$. To make the norm independent of the choice of $\bf i$ one may also define the norm
$$
\|f\|_{\mathscr A(\Omega)}=\sup_{{\bf i}\in\mathbb S^2} \|f\|_{{\mathcal A}(\Omega_{\bf i})}.
$$
In any case, the previous discussion immediately shows that on a given
 bounded axially symmetric s-domain $\Omega$ in $\mathbb{H}$,
the linear space $\mathscr A(\Omega)$  endowed with one of the previous norms is complete.
\\
Note also that $\mathscr A(\Omega)$ can be seen as an inner product linear space, if endowed with the inner product (\ref{scalar}) defined in ${\mathcal A}(\Omega_{\bf i})$, for some fixed ${\bf i}\in\mathbb S^2$. This inner product yields to one of the equivalent norms  $\|\cdot\|_{{\mathcal A}(\Omega_{\bf i})}$.

We now give the following definition  (see \cite{CGS3}) which will be used in Section 7:
%
%
\begin{defin}\label{AC def}
 Let $\Omega$ be an axially symmetric s-domain. We define
$$\mathcal A_c(\Omega):=  P_{\bf  i}(Hol_c(\Omega_{\bf i})) \cap \mathcal L_2(\Omega_{\bf i},\mathbb C({\bf i})),$$
where $\mathcal L_2(\Omega_{\bf i},\mathbb C({\bf i}))$ denotes the space of square integrable functions defined on $\Omega_{\bf i}$ and with values in $\mathbb C(\bf i)$.
\end{defin}

\bigskip

  As explained in \cite{CGLSS}, given any $q\in\Omega_{\bf i}$ the evaluation functional $\phi_q: \mathcal A(\Omega_{\bf i}) \to \mathbb H$, given by
  $$\phi_q[f]:=f(q),\quad \forall f\in \mathcal A(\Omega_{\bf i}),$$
  is a bounded quaternionic right linear functional on $\mathcal A(\Omega_{\bf i}) $ for every ${\bf i} \in \mathbb S^2$.
Thus the Riesz representation theorem for quaternionic right linear Hilbert space shows  the existence of the unique  function $K_q({\bf i}) \in \mathcal A(\Omega_{\bf i}) $ such that
  \begin{equation}\label{RK}
  \phi_q[f_{|\Omega_{\bf i}}]=\langle K_q({\bf i}), f_{|\Omega_{\bf i}}\rangle_{\mathcal A(\Omega_{\bf i}) }.
  \end{equation}
  We set $\mathcal K_{\Omega_{\bf i}}(q,\cdot):= \bar K_q({\bf i})$ and we have:
\begin{defin}[Slice regular Bergman kernel  associated with $\Omega_{\bf i}$]
  The function $$\mathcal K_{\Omega_{\bf i}}(\cdot,\cdot):\Omega_{\bf i}\times \Omega_{\bf i} \to \mathbb H$$ will be called the slice  Bergman kernel associated with $\Omega_{\bf i}$.
\end{defin}
This kernel satisfies the expected properties, see Remark 5 in \cite{CGLSS}.
We can now give the following:
\begin{defin}   {\bf  Slice Bergman kernel  of  the  second  kind   associated with $\Omega$.}
Let  $\mathcal K_{\Omega_{\bf i}}(\cdot,\cdot):\Omega_{\bf i}\times \Omega_{\bf i} \to \mathbb H$
be slice Bergman kernel associated with $\Omega_{\bf i}$,
we will call slice  Bergman kernel   of  the  second  kind      associated with $\Omega$ the function
$$
\mathcal K_{\Omega}: \Omega\times \Omega \to \mathbb{H},
$$
$$
\mathcal K_{\Omega}(x+y{\bf{i}}, r):=\frac{1}{2}(1-{\bf i} {\bf j}) \mathcal K_{\Omega_{\bf j}}(x+y{\bf j},r)+ \frac{1}{2}(1+{\bf i} {\bf j})\mathcal K_{\Omega_{\bf j}}(x+y{\bf j},r).
 $$
\end{defin}

Also this kernel satisfies the expected properties, see \cite{CGLSS}, Proposition 7.
In particular, $\mathcal K_{\Omega}(\cdot,\cdot)$  is a reproducing kernel on $\Omega$: let
 ${\bf{i}}\in \mathbb{S}^2$ then
$$
f(q)=\int_{\Omega_{{\bf i}}}\mathcal{K}_{\Omega}(q,\cdot)\, f\,d\sigma_{{\bf i}},
$$
and the integral does not depend on ${\bf{i}}\in \mathbb{S}^2$. This fact follows from the observation that, see (\ref{RK}):
$$f(q)=\langle \overline{\mathcal{K}(q,\cdot)},f(\cdot)\rangle_{\mathcal A(\Omega_{{\bf i}})}= \langle K_q({\bf i}) ,f(\cdot)\rangle_{\mathcal A(\Omega_{{\bf i}})}. $$

\section{The Bergman kernel $\mathcal K_{\Omega}$ on the unit ball and on the half space}

As it has been pointed out in \cite{K}, even in the classical complex case, the Bergman kernel cannot be computed explicitly on an arbitrary domain, at least not in general. However, if one considers the unit disk
with center at the origin in the complex plane $\mathbb{C}$, then it is well known that the Bergman kernel is given by
\begin{equation}\label{BergmanC}
K(z,\zeta)=\frac{1}{\pi}\displaystyle\frac{1}{(1-z\bar\zeta)^2}.
\end{equation}
We will use this  formula  to  obtain   an explicit expression for the Bergman kernel  of  the  second  kind on the unit ball $\mathbb{B}^4$.
\begin{teo}
The slice regular Bergman kernel    of  the  second  kind    for    the unit ball $\mathbb{B}^4$ is
\begin{equation}\label{Bergman1}
\mathcal{K}_{ \mathbb B^4 }   (q,r)= \frac{1}{\pi} (1-2\bar q\bar r +\bar q^2\bar r^2)(1-2 {\rm Re} [ q] \bar r+|q|^2\bar r^2)^{-2}.
\end{equation}
\end{teo}
Proof. Let us consider any $r\in\mathbb{H}$ and let us choose $z$ such that $z,r$ belong to the same complex plane $\mathbb{C}(\mathbf{j})$.
In this case the Bergman kernel for the unit disc in the complex plane $\mathbb{C}(\mathbf{j})$ is the function $K(z,r)$ defined in (\ref{BergmanC}). The Bergman kernel we are looking for is the slice regular extension of $K(z,r)$, which is unique, to the whole unit ball $\mathbb{B}^4$. The extension is given by the formula (\ref{relationflats}). Thus we have
$$
\mathcal{K}_{ \mathbb B^4 }     (q,r)= \frac{1}{\pi} \Big\{ \frac{1}{2}(1-{\bf i} {\bf j}) \displaystyle\frac{1}{(1-(x+\mathbf{j}y)\bar r)^2} +\frac{1}{2}(1+{\bf i} {\bf j}) \displaystyle\frac{1}{(1-(x-\mathbf{j}y)\bar r)^2}\Big\},
$$
where $q=x+\mathbf{i}y$.
After some standard computations we obtain
$$
\mathcal{K}_{ \mathbb B^4 }     (q,r)
=\frac{1}{\pi}  (1-2\bar q\bar r+\bar q^2\bar r^2)(1-2{\rm Re}(q)\bar r+|q|^2\bar r^2)^{-2}.
$$
\hfill $\blacksquare$

The regularity properties  of the Bergman kernel of the second kind on the unit ball can be proved directly:
\begin{pro}
The kernel $\mathcal{K}_{ \mathbb B^4 }     (q,r)$   is slice regular in $q$ and slice right anti-regular in the variable $r$.
\end{pro}
Proof.
By construction, the Bergman kernel is slice regular in the variable $q$. The fact that $\mathcal{K}_{ \mathbb B^4 }   (q,r)$ is slice right anti-regular in the variable $r$ can be verified directly: the function $p(q,r):=(1-2\bar q\bar r +\bar q^2\bar r^2)$ is  a polynomial in the variable $\bar r$ with coefficients on the left thus it is slice right anti-regular in $r$; the function $h(q,r):=(1-2 {\rm Re} ( q) \bar r+|q|^2\bar r^2)^{-2}$ is   a   rational function in $r$ with real coefficients thus it is slice right anti-regular in $r$. Then $\mathcal{K}(q,r)=p(q,r)h(q,r)$ is the product of two slice right anti-regular functions in $r$, where $h$, as a function in the variable $r$, has real coefficients. The statement follows.

\hfill $\blacksquare$

\bigskip

In the paper \cite{csTrends}, the authors show that the slice hyperholomorphic Cauchy kernel admits two different analytic expressions (see also \cite{CGS} for the specific case of regular functions). An interesting feature of the slice regular Bergman kernel for   the unit ball is that it also admits two analytic expressions. In fact, we have the following result.
\begin{pro}
The kernel $\mathcal{K}_{ \mathbb B^4 }   (q,r)$   can be written in the form
\begin{equation}\label{Bergman2}
\mathcal{K}_{ \mathbb B^4 }    (q,r)=\frac{1}{\pi}(1-2q{\rm Re}[r]+q^2|r|^2)^{-2}(1-2qr+q^2r^2).
\end{equation}
\end{pro}

Proof.
First of all, we note that the function
$$
\psi(q,r):=\frac{1}{\pi}(1-2q{\rm Re}[r]+q^2|r|^2)^{-2}(1-2qr+q^2r^2)
$$
is slice regular in the variable $q$ since it is the product of a rational function with real coefficients and a polynomial in $q$ with coefficients on the right. Then we observe that for any fixed $r=u+v\mathbf{i}$ if we consider $q$ belonging to the plane $\mathbb{C}(\mathbf{i})$ we have:
$$
\mathcal{K}_{ \mathbb B^4 }   (q,r)_{|\mathbb{C}(\mathbf{i})}=\psi(q,r)_{|\mathbb{C}(\mathbf{i})}=\displaystyle\frac{1}{(1-q\bar r)^2}.
$$
By the Identity Principle, the two functions $\mathcal{K}_{ \mathbb B^4 } $ and $\psi$ coincide (and $\psi$ turns out to be slice right anti-regular in $r$).
\hfill $\blacksquare$

\bigskip

Let us now consider the half space $\mathbb{H}^+=\{q\in\mathbb{H}\ |\ {\rm Re}(q)>0 \}$. To construct the Bergman kernel associated with $\mathbb H^+$,
  we work in the complex plane and note that the transformation $\alpha(z)=\frac{z-1}{z+1}$ maps the complex half plane of numbers with positive real part  $\mathbb{C}_+$ to the unit disc $\mathbb{D}$. By a well known formula, see \cite{vasilevski}, p. 37, we have
that the Bergman kernel of $\mathbb{C}_+$ is the function ${K}_{\mathbb{C}_+}=\frac{1}{\pi}\frac{1}{(z+\bar\xi)^2}$. We will extend this function in order to obtain
the Bergman kernel of $\mathbb{H}^+$.

\begin{teo}
The slice regular Bergman kernel    of  the  second  kind    for    the  half space $ \mathbb{H}^+$ is
\begin{equation}\label{Bergmanhalf}
\begin{split}
\mathcal{K}_{ \mathbb{H}^+}   (q,r) &= \frac{1}{\pi} (\bar{q}^2+2\bar q\bar r +\bar r^2)(|q|^2 +2 {\rm Re} [ q] \bar r+\bar r^2)^{-2}\\
&=\frac{1}{\pi} (q^2 +2 {\rm Re} [ r] q +| r|^2)^{-2}({q}^2+2 q r + r^2).\\
\end{split}
\end{equation}
\end{teo}
Proof. The proof is based on computations similar to those done to obtain formulas (\ref{Bergman1}) and (\ref{Bergman2}), thus we will not repeat them.\hfill $\blacksquare$

\begin{ob}{\rm
Another possibility to prove the result, but with more complicated computations, is to use (\ref{Bergman1}) and Proposition 4.2 in \cite{CGS2}.
}
\end{ob}

\section{Comparison between the slice regular Bergman spaces of first and second kind}
The  properties presented so far are satisfied  in a domain  $\Omega\subset \mathbb{H}$, however   the concept of slice regular function deals with the fact that the restriction of a function to  open sets $\Omega_{\mathbf{i}}\subset \mathbb{C}(\mathbf{i})$ is a holomorphic map. Therefore the behavior of a slice regular function is related with that one of its restrictions to the sets  $\Omega_{\mathbf{i}}$.
Purpose of this section is to compare the sets $\mathcal A(\Omega)$ and $\mathcal A(\Omega_{\bf i})$. We will show that they contain different elements unless a suitable weight function is taken into account.
\begin{defin}
Let $\Omega$ be a bounded axially symmetric s-domain.
\begin{enumerate}
\item For any  ${\bf i}\in \mathbb S^2$ define the weight function  $\rho(z)=Im(z)^2,\quad \forall z\in \Omega_{\bf i}$ and the weighted Bergman space on $\Omega_{\bf i}$  with weight $\rho$: $$ \mathcal A_{\rho}(\Omega_{\bf i})  := \{  f\in \mathcal{SR}(\Omega) \ \mid \ \int_{\Omega_{\bf i}} |f_{\mid_{\Omega_{\bf i} }} |^2  \rho d\sigma_{\bf i} <\infty  \}  ,$$
endowed with the inner product
$$\langle f, g \rangle_{\mathcal A_{\rho}(\Omega_{\bf i}) }:=  \int_{\Omega_{\bf i}} \bar f g  \rho d\sigma_{\bf i},$$
and the corresponding norm
  $$\|f\|_{\mathcal A_{\rho}(\Omega_{\bf i}) }:=  \left( \int_{\Omega_{\bf i}} |f_{\mid_{\Omega_{\bf i} }}|^2  \rho d\sigma_{\bf i}\right)^2.$$
\item Define the weight function  $\delta(q)= \frac{1}{|\underline{q}|^2}$, and the weighted Bergman space on $\Omega$ with weight $\delta$: $$ \mathcal A_{\delta}(\Omega)  := \{  f\in \mathcal{SR}(\Omega) \ \mid \ \int_{\Omega} |f |^2  \delta d\mu <\infty  \}  ,$$
endowed with the inner product
$$\langle f, g \rangle_{\mathcal A_{\delta}(\Omega) }:=  \int_{\Omega} \bar f g  \delta d\mu,$$
and the corresponding norm $$\|f\|_{\mathcal A_{\delta}(\Omega) }:=  \left( \int_{\Omega} |f|^2  \delta d\mu\right)^2.$$
\end{enumerate}
\end{defin}
\begin{pro}\label{relationdomain-slice}  Let $f\in \mathcal A(\Omega)$ and ${\bf i}$ be a fixed element in $\mathbb S^2$. Then
\begin{equation}\label{inequality} \frac 14\int_{\Omega}|f(x+{\bf i}y)|^2 d\mu \leq  \|f\|_{\Omega}^2 \leq   4 \int_{\Omega}|f(x+{\bf i}y)|^2 d\mu .\end{equation}
\end{pro}
Proof. To prove the result, we use the fact that for any ${\bf i}, {\bf j} \in\mathbb S^2$ one has that $|f(x+{\bf i}y)|^2 \leq 2[ |f(x+{\bf j}y)|^2+|f(x-{\bf j}y)|^2]$, see (\ref{inequalitymodules}). Let us take the integral on $\Omega$ of both sides of this inequality and let us consider the right hand side.
A change of variable in the second integral below and the axial symmetry of $\Omega$ give:
\begin{equation}\label{ineint}
\int_\Omega  |f(x+{\bf j} y)|^2 d\mu +\int_\Omega  |f(x-{\bf j} y)|^2 d\mu=2\int_\Omega  |f(x+{\bf j} y)|^2 d\mu.\\
\end{equation}
Let $q=x+{\bf j} y$ vary in $\Omega$, and keep ${\bf i}\in\mathbb S^2$ fixed. Then (\ref{ineint}) and (\ref{inequalitymodules}) give
$$
 \frac 14 \int_\Omega |f(x+{\bf i} y)|^2 d\mu \leq \int_\Omega  |f(x+{\bf j} y)|^2 d\mu  = \|f\|_{\Omega}^2 .
$$
Exchanging the role of ${\bf i}$ and ${\bf j}$ and repeating the reasoning, we obtain
$$
 \|f\|_{\Omega}^2 = \int_\Omega |f(x+ {\bf j} y)|^2 d\mu \leq 4\int_\Omega  |f(x+ {\bf i} y)|^2 d\mu ,
$$
and the statement follows.
\hfill $\blacksquare$
\\
\begin{pro}  Let $\bf{i} $ be a fixed element in $\mathbb S^2$.
 \begin{enumerate}
 \item If  $f \in \mathcal A(\Omega)$ then there  exist   $K_\Omega>0$
 such that
$$ \int_{\Omega_{\bf i}}|f_{\mid_{\Omega_{\bf i}}}|^2 \rho d\sigma_{\bf i} \leq K_\Omega   \|f\|_{\mathcal  A(\Omega) }^2   ;$$
if $f\in \mathcal A_{\rho}(\Omega_{\bf i})$ then there exists      $M_{\Omega} > 0$ such that
$$ \|f\|_{\mathcal  A(\Omega) }^2
\leq M_\Omega   \int_{\Omega_{\bf i}}|f_{\mid_{\Omega_{\bf i}}}|^2 \rho d\sigma_{\bf i}.$$
 \item  If  $f \in \mathcal A(\Omega_{\bf i})$ then there  exist  $ {K'}_\Omega> 0$  such that
$$ \int_{\Omega}|f|^2 \delta d\mu  \leq {K'}_{\Omega}  \|f\|_{\mathcal  A(\Omega_{\bf i}) }^2;$$
if $f\in \mathcal A_{\delta }(\Omega) $ then there exists    ${M'}_{\Omega}> 0$ such that
$$ \|f\|_{\mathcal  A(\Omega_{\bf i}) }^2  \leq    {M}_\Omega'\int_{\Omega}|f|^2 \delta d\mu .$$
  \end{enumerate}\end{pro}

  Proof.
\begin{enumerate}
\item  Assume that $f \in \mathcal A(\Omega)$.  Consider ${\bf i}\in\mathbb S^2$ fixed and the integral
$$
\int_\Omega  |f (q_0 +  |\underline{q}| {\bf i}) |^2 d\mu .
$$
This integral exists, finite, by virtue of our assumption and by (\ref{inequality}).
Consider the change of coordinates $q=q_0+ {\bf i}\, |\underline{q}| = r\cos\theta + {\bf i} \ r \sin\theta$, where ${\bf i}=(\sin\delta\cos\varphi, \sin\delta\sin\varphi, \cos\delta)$, and the values of $r,\theta,\delta, \varphi$ are such that $q\in\Omega$. Note that the Jacobian of the change of coordinate is $r^3\sin^2\theta \sin\delta$. We have that
\begin{equation}\label{equiv1}
\int_\Omega  |f (q_0 +  |\underline{q}| {\bf i}) |^2 d\mu   =      \lambda_\Omega   \int_r \int_\theta  | f_{\mid_{\Omega_{\bf i}}} |^2  (r\sin\theta)^2 rdr\ d\theta=  \lambda_\Omega   \int_{\Omega_{\bf i}}  | f_{\mid_{\Omega_{\bf i}}} |^2  \rho  d\sigma_{\bf i} .
\end{equation}
 The first part of the statement is a direct consequence of  (\ref{equiv1}) and Proposition \ref{relationdomain-slice}.
 To get the second part, assume that $f\in \mathcal A_{\rho}(\Omega_{\bf i})$. Using backward the equalities in (\ref{equiv1}) together with Proposition \ref{relationdomain-slice} the statement follows.
\item Let $f\in \mathcal A (\Omega_{\bf i})$ with ${\bf i}\in\mathbb S^2$ fixed. Let us compute
$$\int_\Omega  |f (q_0 +  |\underline{q}| {\bf i}) |^2\left(\frac{1}{y}\right)^2 d\mu . $$
Using the change of coordinate described in point 1 we have:
 $$\int_\Omega  |f (q_0 +  |\underline{q}| {\bf i}) |^2\left(\frac{1}{y}\right)^2 d\mu  =      \lambda_\Omega   \int_r \int_\theta  | f_{\mid_{\Omega_{\bf i}}} |^2\left(\frac{1}{r\sin\theta}\right)^2  (r\sin\theta)^2 rdr\ d\theta=  \lambda_\Omega   \int_{\Omega_{\bf i}}  | f_{\mid_{\Omega_{\bf i}}} |^2   d\sigma_{\bf i} , $$
 where $\lambda_{\Omega}$ is a constant depending of $\Omega$.
 Reasoning as in the proof of Proposition \ref{relationdomain-slice} we deduce
 $$\int_\Omega|f|^2\delta d\mu  \leq   4\int_{\Omega}|f(x+{\bf i}y)|^2\left(\frac{1}{y}\right)^2  d\mu, $$
from which the first part of the statement follows.
\\
On the other hand, consider $f\in \mathcal A_{\delta}(\Omega)$, then
$$ \int_{\Omega}|f(x+{\bf i}y)|^2\left(\frac{1}{y}\right)^2 d\mu_q \leq  4\int_\Omega|f|^2\delta d\mu  .$$
Using the previous inequalities one has the result.
\end{enumerate}
\hfill $\blacksquare$

\begin{cor}
      The sets of functions $\mathcal A_{\rho}(\Omega_{\bf i})$ and $ \mathcal A(\Omega)$ contain the same elements and have equivalent norms. Similarly, for the sets  $ \mathcal A(\Omega_{\bf i})$ and $ \mathcal A_\delta(\Omega)$.
      \end{cor}
\begin{ob}{\rm
Consider the unit ball  $\mathbb B^4:=\{ q\in \mathbb H \ \mid \ |q|<1  \}$, the unit disc $\mathbb D_{\bf i}=\mathbb B\cap\mathbb C({\bf i})$ and the following slice regular function on $\mathbb B^4$:
$$f(q) = \frac{1}{1-q}, \quad \forall q\in \mathbb B^4.$$
For any ${\bf i} \in \mathbb S^2$  the integral
$$\int_{\mathbb D_{\bf i}} \frac{1}{|1-z |^2}d\sigma_{\bf i} $$
is not finite. However, both the integrals
$$\int_{\mathbb D_{\bf i}} \frac{1}{|1-z |^2} ( y^2 )d\sigma_{\bf i},$$
$$\int_{\mathbb B^4} \frac{1}{|1-q |^2}d\mu $$
are finite.
Therefore $f\in \mathcal A(\mathbb B^4) \setminus \mathcal A(\mathbb D_{\bf i})$ and $f\in\mathcal A_\rho (\mathbb D_{\bf i})$.
}
\end{ob}

\section{The Bergman-Fueter transform on the unit ball}

In this section we introduce an integral transform which associates to every slice regular function $f$  defined on the unit ball $\mathbb{B}^4$ of $\mathbb{H}$ a Fueter regular function
$\breve{f}$ defined on the same set. We will call  the   mapping   $f\to \breve{f}$ slice Bergman-Fueter integral transform
(for short BF-integral transform). Its definition is based on the Fueter mapping theorem and on the slice regular Bergman kernel  of  the  second kind.
It is inspired by the paper \cite{CoSaSo}, where
 an integral transform that generates Fueter regular functions from slice regular functions is defined. Precisely, in \cite{CoSaSo} we
consider the function
$$
\mathcal{F}(s,q):=\Delta S^{-1}(s,q)=-4(s-\bar q)(s^2-2{\rm Re}[q]s +|q|^2)^{-2},
$$
where    $\Delta$    is  the  Laplace  operator  and   $S^{-1}(s,q)$ is the slice regular Cauchy kernel.
Let $W\subset \mathbb{H}$ be an axially symmetric open set and let $f$ be  a slice regular function on $W$. Let $\Omega$ be a
bounded  axially symmetric open set such that $\overline{\Omega}\subset W$. Suppose that the boundary of $\Omega\cap \mathbb{C}(\mathbf{i})$  consists of a
finite number of continuously differentiable Jordan curves for any $\mathbf{i}\in\mathbb{S}^2$.
Then, if $q\in \Omega$, the Cauchy-Fueter regular function $\breve{f}(q)$     given by
$$
\breve{f}(q):=\Delta f(q)
$$
has the integral representation
\begin{equation}\label{Fueter_quat}
 \breve{f}(q)=\frac{1}{2 \pi}\int_{\partial(\Omega\cap\mathbb{C}(\mathbf{i}))} \mathcal{F}(s,q)ds_\mathbf{i} f(s),\ \ \ ds_\mathbf{i}=-ds  \mathbf{i},
\end{equation}
 and the
integral does not depend  neither   on $\Omega$ and  nor   on the imaginary unit
$\mathbf{i}\in\mathbb{S}^2$.

In the case under study, we restrict our attention to the case of the unit ball $\mathbb B^4$,  the   case in which we know the explicit
expression of the Bergman kernel, but such definition can be given for more general domains. We make the following preliminary observations.
\begin{itemize}
\item
 It is well known that  slice regular functions on axially symmetric s-domains are of class $C^{\infty}$.
\item
The Fueter mapping theorem, for $f\in \mathcal A(\mathbb B^4)$ gives:
$$
    \breve{f}(q)=\int_{\partial(\mathbb B^4\cap\mathbb{C}(\mathbf{i}))} \Delta \mathcal{K}_{\mathbb B^4}  (q,r) f(r) d\sigma(r).
  $$
\end{itemize}

\begin{defin}
The integral
   \begin{equation}\label{intrsB}
    \breve{f}(q)=\int_{\partial(\mathbb B^4\cap\mathbb{C}(\mathbf{i}))} \Delta \mathcal{K}_{\mathbb B^4}  (q,r) f(r) d\sigma(r).
  \end{equation}
   is  called  the   Bergman-Fueter  integral transform of $f$ and  its  kernel
  $$
  \mathcal{K}_{BF}(q,r):=\Delta \mathcal{K}_{\mathbb B^4}  (q,r)
  $$
  is called Bergman-Fueter kernel;   here
$\mathcal{K}_{\mathbb B^4}   (q,r)$ is the slice  Bergman kernel of  the  second  kind  for   the unit ball $\mathbb{B}^4$.
\end{defin}
 We will study in another paper the Sobolev regularity of the Bergman-Fueter integral transform.

\begin{teo}
The following formula holds:
$$
\Delta \mathcal{K}_{\mathbb B^4}  (q,r)=-\frac{4}{\pi}[(1-2 Re[q] \overline{r})+|q|^2\overline{r}^2)^{-2}+2(1-2\bar{q}\,\bar{r}+
\overline{q}^2\overline{r}^2)(1-2 Re[q] \overline{r})+|q|^2\overline{r}^2)^{-3}] \overline{r}^2 .
$$
 \end{teo}
Proof. To compute the derivatives of the Bergman kernel written in the form (\ref{Bergman1}), we set for simplicity
$$
\mathcal{\widetilde{K}}(q,r):=\pi\mathcal{{K}}_{\mathbb B^4}   (q,r)=(1-2\bar q\bar r +\bar q^2\bar r^2)(1-2 {\rm Re} [ q] \bar r+|q|^2\bar r^2)^{-2}.
$$
So we have:
$$
\pp_{x_0}\mathcal{\widetilde{K}}(q,r)=
(-2\overline{r}+2\overline{q}\, \overline{r}^2)(1-2 Re[q] \overline{r})+|q|^2\overline{r}^2)^{-2}
$$
$$
-2(1-2\bar q\bar r +\bar q^2\bar r^2)(-2\overline{r}+2x_0 \overline{r}^2)(1-2 {\rm Re} [ q] \bar r+|q|^2\bar r^2)^{-3}
$$
and
$$
\pp^2_{x_0}\mathcal{\widetilde{K}}(q,r)=2 \overline{r}^2(1-2 Re[q] \overline{r})+|q|^2\overline{r}^2)^{-2}
$$
$$
-4(-2\overline{r}+2\overline{q}\, \overline{r}^2)(-2\overline{r}+2x_0 \overline{r}^2)(1-2 {\rm Re} [ q] \bar r+|q|^2\bar r^2)^{-3}
$$
$$
-4(1-2\bar q\bar r +\bar q^2\bar r^2)\overline{r}^2(1-2 {\rm Re} [ q] \bar r+|q|^2\bar r^2)^{-3}
$$
$$
+6(1-2\bar q\bar r +\bar q^2\bar r^2)(-2\overline{r}+2x_0 \overline{r}^2)^2(1-2 {\rm Re} [ q] \bar r+|q|^2\bar r^2)^{-4}.
$$
We now calculate the first and the second derivatives with respect to $x_1$ We have:
$$
\pp_{x_1}\mathcal{\widetilde{K}}(q,r)=
(2 {\imath}\bar r-{\imath}\bar q\bar r^2-\bar q {\imath}\bar r^2)(1-2 {\rm Re} [ q] \bar r+|q|^2\bar r^2)^{-2}
$$
$$
-2(1-2\bar q\bar r +\bar q^2\bar r^2) \, 2x_1 \overline{r}^2  \,(1-2 {\rm Re} [ q] \bar r+|q|^2\bar r^2)^{-3}.
$$
The second derivative with respect to $x_1$  is
$$
\pp^2_{x_1}\mathcal{\widetilde{K}}(q,r)=
-2  \overline{r}^2(1-2 {\rm Re} [ q] \bar r+|q|^2\bar r^2)^{-2}
$$
$$
-8(2 {\imath}\bar r-{\imath}\bar q\bar r^2-\bar q {\imath}\bar r^2)  \, x_1 \overline{r}^2\, (1-2 {\rm Re} [ q] \bar r+|q|^2\bar r^2)^{-3}
$$
$$
-4(1-2\bar q\bar r +\bar q^2\bar r^2) \,  \overline{r}^2  \,(1-2 {\rm Re} [ q] \bar r+|q|^2\bar r^2)^{-3}
$$
$$
+24(1-2\bar q\bar r +\bar q^2\bar r^2) \, x_1^2 \overline{r}^4  \,(1-2 {\rm Re} [ q] \bar r+|q|^2\bar r^2)^{-4}.
$$
The derivatives with respect to
 $x_2$ and $x_3$ can be computed in a similar way.
With some calculations we obtain:
$$
\Delta \mathcal{\widetilde{K}}(q,r)
=-4\overline{r}^2(1-2 {\rm Re} [ q] \bar r+|q|^2\bar r^2)^{-2}
$$
$$
-16\Big( 1 -  2 \overline{q}\, \overline{r}+\overline{q}^2\, \overline{r}^2\Big)\overline{r}^2(1-2 {\rm Re} [ q] \bar r+|q|^2\bar r^2)^{-3}
-16(1-2\bar q\bar r +\bar q^2\bar r^2)\overline{r}^2(1-2 {\rm Re} [ q] \bar r+|q|^2\bar r^2)^{-3}
$$
$$
+24(1-2\bar q\bar r +\bar q^2\bar r^2)\Big(
1-2{\rm Re} [ q]\overline{r} +|q|^2\overline{r}^2\Big) \overline{r}^2(1-2 {\rm Re} [ q] \bar r+|q|^2\bar r^2)^{-4}.
$$
We finally get:
$$
\Delta \mathcal{\widetilde{K}}(q,r)
=-4\overline{r}^2(1-2 {\rm Re} [ q] \bar r+|q|^2\bar r^2)^{-2}
-8( 1 -  2 \overline{q}\, \overline{r}+\overline{q}^2\, \overline{r}^2)\overline{r}^2(1-2 {\rm Re} [ q] \bar r+|q|^2\bar r^2)^{-3}
$$
and  the statement follows.
\hfill$\blacksquare$

\bigskip
We rewrite the kernel $\Delta \mathcal{\widetilde{K}}(q,r)$ in a shorter way defining the function
$$
\mathcal{Q}(q,r):=(1-2 Re[q] \overline{r})+|q|^2\overline{r}^2)^{-1}.
$$
\begin{cor}\label{strudelberg}
Let $\mathcal{Q}(q,r)$ be as above. Then the kernel $\Delta \mathcal{\widetilde{K}}(q,r)$ is given by
$$
\Delta \mathcal{K}_{\mathbb B^4}  (q,r)=-\frac{4}{\pi}[\mathcal{Q}(q,r)+2\mathcal{K}_{\mathbb B^4}  (q,r)] \mathcal{Q}(q,r)\overline{r}^2.
$$
 \end{cor}   \hfill $\blacksquare$

\bigskip

\begin{teo}
The Bergman-Fueter kernel $\mathcal{K}_{BF}(q,r)$ is Fueter regular in $q$ and slice anti-regular in $r$.
\end{teo}
Proof.
The fact that $\mathcal{K}_{BF}(q,r)$ is Fueter regular in $q$ is a direct consequence of the Fueter mapping theorem since
by definition it is $\mathcal{K}_{BF}(q,r):=\Delta \mathcal{K}_{ \mathbb B^4}  (q,r)$   and  $\mathcal{K}(q,r)$ is slice regular in the variable $q$. We have to verify that it is
 slice anti-regular in $r$. In fact observe that the function $1-2 Re[q] \overline{r})+|q|^2\overline{r}^2$    is  slice anti-regular in $r$ and it has real coefficients so also
 $\mathcal{Q}(q,r):=(1-2 Re[q] \overline{r})+|q|^2\overline{r}^2)^{-1}$ is  slice anti-regular in $r$ and for the same reason also
 $\mathcal{Q}^2(q,r)$. Since
$\mathcal{K}_{\mathbb B^4}  (q,r)$ is slice anti-regular in $r$ then the product  $\mathcal{K}(q,r)\mathcal{Q}(q,r)$ is  slice anti-regular in $r$
because $\mathcal{Q}(q,r)$ is a rational function of a polynomial with real coefficients.
The statement follows from Corollary \ref{strudelberg}.

\hfill $\blacksquare$

\section{Schwarz reflection principle}
In the theory of one complex variable the Schwarz reflection principle is well known. It shows how to extend functions defined on domains in the upper half plane and with real boundary values on the real axis to domains  which  are symmetric with respect to
the real axis. Here we use this property to obtain results for slice regular functions and the Bergman theory.

\noindent
We begin by recalling some fact in classical complex analysis. Consider a domain $\Omega\subset \mathbb C$  such that $\Omega\cap \mathbb R \neq \emptyset $ and for any $z\in \Omega$ one has that  $\bar z \in \Omega$, or equivalently, $\overline\Omega=\Omega$.  Set
 $$\mathbb C^+:=\{z\in \mathbb C \ \mid \ {\rm Im}( z) >0\},\ {\rm and}\ \mathbb C^-:=\{z\in \mathbb C \ \mid \ {\rm Im}( z) < 0\},$$  and $\Omega^+ := \Omega\cap \mathbb C^+$, $\Omega^- := \Omega\cap \mathbb C^-  $.

Recall that the spaces $Hol(\Omega)$ and $Hol_c(\Omega)$ have been defined in Section 2.

\begin{defin}  \label{def2} Let $\Omega\subset\mathbb C$ be such that $\overline\Omega=\Omega$. We define:
  \begin{enumerate}\item  $\displaystyle \widetilde{Hol} (\Omega^+)= \{  f\in Hol(\Omega^+) \ \mid \   \textrm{   for any  } \  x\in\mathbb R  \  \textrm{there exists}   \lim_{\Omega^+\ni z \to x } f(z) \in\mathbb R   \} ;$
\item   the function $$E[f](z)=\left\{ \begin{array}{ll}  f(z) &  \textrm{ if  }  z\in \Omega^+\\ \\
\displaystyle  \lim_{\Omega^+\ni w \to z } f(w) \in\mathbb R   & \textrm{ if } z\in \Omega\cap \mathbb R
\\ \\
\overline{f(\bar z)} & \textrm{ if } z\in \Omega^-  \end{array}  \right.$$
for any $f\in \widetilde{Hol} (\Omega^+)$ ;
\item  $E( \ \widetilde{Hol}(\Omega^+)\  )=\{ E[f] \ \mid \ f\in \widetilde{Hol} (\Omega^+) \}.$\end{enumerate}
\end{defin}

\begin{ob}\label{Sprinciple}{\rm
With the notations in Definition \ref{def2}, the Schwarz reflection principle asserts that
$Hol_c(\Omega) = E( \ \widetilde{Hol}(\Omega^+)\  )$.

\noindent
Denoting by ${\bf i}$ the imaginary unit of the complex plane $\mathbb C$, we have that $Hol(\Omega)= Hol_c(\Omega)  + Hol_c(\Omega){\bf i}$, and so we immediately obtain the formula
$$Hol(\Omega)= E( \ \widetilde{Hol}(\Omega^+)\  ) + E( \ \widetilde{Hol}(\Omega^+)\  ){\bf i}.$$
The previous formula says that any holomorphic function $f\in Hol(\Omega)$ can be written in terms of two elements of $f_1,f_2\in \widetilde{Hol} (\Omega^+)$ as  $ f = E(f_1) + E(f_2){\bf i},$
where
$$f_1(z)=\frac{1}{2} \left( f\mid_{\Omega^+}(z)+\overline{f\mid_{\Omega^+}(\bar z)} \right) ,\quad \textrm{and} \quad
f_2(z)=\frac{-{\bf i}}{2} \left( f\mid_{\Omega^+}(z)-\overline{f\mid_{\Omega^+}(\bar z)} \right) , \quad \forall z\in \Omega^+.$$
}
\end{ob}
Let us now consider the function spaces
\[
\begin{split}
&\widetilde{\mathcal A}(\Omega^+) =\widetilde{Hol} (\Omega^+) \cap \mathcal L_2(\Omega^+) \\
&\mathcal A(\Omega) =Hol (\Omega) \cap \mathcal L_2(\Omega)\\
&\mathcal A_c(\Omega) =Hol_c (\Omega) \cap \mathcal L_2(\Omega)\\
\end{split}
\]
Obviously, $\mathcal A_c(\Omega) \subset \mathcal A(\Omega) $.
The facts
 $Hol(\Omega)= Hol_c(\Omega)  + Hol_c(\Omega){\bf i}$ and   $Hol_c(\Omega) = E( \ \widetilde{Hol}(\Omega^+)\  )$ imply that
 $$\mathcal A(\Omega) = E( \ \widetilde{\mathcal A}(\Omega^+)  \  ) + E( \ \widetilde{\mathcal A}(\Omega^+) \  ){\bf i}.$$
The following results are presented with their proofs for the sake of completeness and for lack of reference.
\begin{pro}  \label{pro5}    {   }  Let $\Omega\subseteq\mathbb C$ be an open set such that $\bar\Omega=\Omega$. Then:
\begin{enumerate}
\item  If $f,g\in \mathcal A_c(\Omega)$   then
  $\displaystyle  <f, g>_{\mathcal A(\Omega)} =2 \ {\rm Re}  \   \left( \int_{\Omega^+} \bar  f  g       d\sigma \right)$.
\item Let $f,g \in \mathcal A(\Omega)$   and let  $f_1,f_2,g_1,g_2 \in \mathcal A_c(\Omega) $  be such that  $f=f_1 + f_2{\bf i}$  and $g=g_1 +g_2{\bf i}  $     then
  $$\displaystyle  <f, g>_{\mathcal A(\Omega)} =  2\left[   {\rm Re}   \ \left(\int_{\Omega^+} (\bar f_1 g_1 +\bar  f_2 g_2 ) d\sigma \right)  + {Re}\ \left( \int_{\Omega^+}  ( \bar f_1 g_2  -  \bar f_2 g_1)d\sigma       \right)  {\bf i}\right]$$
  \item If $f\in \mathcal A_c(\Omega)$, then $$\|f\|_{\mathcal A(\Omega)} =  \sqrt{2}  \left[  \int_{\Omega^+} |f|^2 d\sigma\right]^{\frac{1}{2}},$$
and if  $f \in \mathcal A(\Omega)$ is written as  $f=f_1 + f_2{\bf i}$ with  $f_1,f_2\in \mathcal A_c(\Omega) $, then
  $$\|f\|_{\mathcal A(\Omega)} =  \sqrt{2}  \left[  \int_{\Omega^+} ( \ |f_1|^2 + |f_2|^2 \ ) d\sigma\right]^{\frac{1}{2}}. $$
\end{enumerate}
     \end{pro}
Proof.
 To prove point 1  we compute
 $$\displaystyle  <f, g>_{\mathcal A(\Omega)} =  \int_{\Omega} \bar  f  g       d\sigma =\int_{\Omega^+} \bar  f  g       d\sigma +\int_{\Omega^-} \bar  f  g       d\sigma$$
$$=\int_{\Omega^+} \bar  f(\zeta)  g(\zeta)       d\sigma+ \int_{\Omega^-}   f(\bar \zeta)  \overline{ g(\bar \zeta)}        d\sigma=
\int_{\Omega^+} \bar  f (\zeta)  g (\zeta)      d\sigma+ \int_{\Omega^+}  \overline{ \overline{f( \zeta)}   g( \zeta)}        d\sigma=
 2 {\rm Re}  \left( \int_{\Omega^+} \bar  f  g       d\sigma \right).$$
Points 2. and 3. are consequence of 1.
\hfill $\blacksquare$

In the next result $\mathcal B_{\Omega}(\cdot, \cdot )$ denotes the Bergman kernel associated with $\Omega$.
By Remark \ref{Sprinciple}, there exist
$\mathcal B_{\Omega,1}, \ \mathcal B_{\Omega,2}:\Omega\times \Omega^+\to \mathbb C$ be such that
$$\mathcal B_{\Omega}(\cdot, z  ) =    E(\mathcal B_{\Omega,1}( \cdot,z) ) + E(\mathcal B_{\Omega,2}(\cdot, z )){\bf i} .$$
\begin{pro}  Let $\Omega\subseteq\mathbb C$ be an open set such that $\Omega=\overline{\Omega}$, let $z\in \Omega\setminus \Omega_0$, and
let  $\mathcal B_{\Omega,1}, \ \mathcal B_{\Omega,2}:\Omega\times \Omega^+\to \mathbb C$ be such that
$$\mathcal B_{\Omega}(\cdot, z  ) =    E(\mathcal B_{\Omega,1}( \cdot,z) ) + E(\mathcal B_{\Omega,2}(\cdot, z )){\bf i} .$$
\begin{enumerate}
\item  For any   $f \in \mathcal A(\Omega)$, let  $f_1,f_2\in \mathcal A_c(\Omega) $ be  such that  $f=f_1 + f_2{\bf i}$. Then

$$\begin{array}{l}\displaystyle f(z)=2\left[ {\rm  Re}  \  \left( \int_{\Omega^+}( \ \mathcal B_{\Omega,1}(z,\zeta)f_1(\zeta) + \mathcal B_{\Omega,2}(z,\zeta)f_2(\zeta) \  )d\sigma   \right)   \right.
+   \\
\\
 \  \hspace{5cm}  \left.\displaystyle +  {\rm Re} \ \left( \int_{\Omega^+}(  \  \mathcal B_{\Omega,1}(z,\zeta)f_2(\zeta) - \mathcal B_{\Omega,2}(z,\zeta)f_1(\zeta)  \ )d\sigma   \right)  \ \bf{i}.
\right]\end{array}$$

\item If  $u,v$ are the real component functions of $f$, that is  $f=u+ v\bf{i}$, then
$$\begin{array}{l}\displaystyle u(z)=2 {\rm Re}  \  \left( \int_{\Omega^+}( \ \mathcal B_{\Omega,1}(z,\zeta)f_1(\zeta) + \mathcal B_{\Omega,2}(z,\zeta)f_2(\zeta) \  )d\sigma   \right)
  \\
\\
v(z) = 2 \displaystyle   {\rm Re} \ \left( \int_{\Omega^+}(  \  \mathcal B_{\Omega,1}(z,\zeta)f_2(\zeta) - \mathcal B_{\Omega,2}(z,\zeta)f_1(\zeta)  \ )d\sigma   \right) .  \end{array}$$
\item If $f\in \mathcal A_c(\Omega)$ then
$$f(z) = 2\left[ {\rm  Re} \left( \int_{\Omega^+}     \mathcal B_{\Omega, 1}(z,\zeta) f(\zeta) d\sigma  \right)  - {\bf{i}}  \int_{\Omega^+}  \mathcal B_{\Omega,2}(z,\zeta)   f(\zeta)  d\sigma \right].
  $$
\end{enumerate}
\end{pro}
Proof.
\begin{enumerate}
\item It is a consequence of Proposition \ref{pro5}.
\item Point 2 follows from  1. by direct computations.
\item
$$\begin{array} {rl}f(z) =&\displaystyle  \int_{\Omega^+}     \left[ \left( \mathcal B_{\Omega,1}(z,\zeta)-  \mathcal B_{\Omega,2}(z,\zeta) {\bf{i}} \right) f(\zeta)   +
\left ( \overline{ \mathcal B_{\Omega,1}(z,\zeta)    + {\bf{i}} \mathcal B_{\Omega,2}(z,\zeta)} \right) \bar f(\zeta) \right] d\sigma
  \\
&\\
=&\displaystyle  \int_{\Omega^+}     \left[ \left( \mathcal B_{\Omega,1}(z,\zeta)+  \mathcal B_{\Omega,2}(z,\zeta) {\bf{i}} \right) f(\zeta)     -2  \mathcal B_{\Omega,2}(z,\zeta){\bf{i}}  f(\zeta)     +
 \overline{ \mathcal B_{\Omega}(z,\zeta)  f(\zeta) }\right] d\sigma\\
&\\
=&\displaystyle \int_{\Omega^+}     \left[ \mathcal B_{\Omega}(z,\zeta) f(\zeta)     -2  \mathcal B_{\Omega,2}(z,\zeta) {\bf{i}}  f(\zeta)     +
 \overline{ \mathcal B_{\Omega}(z,\zeta)  f(\zeta) }\right] d\sigma
  \end{array}$$
\end{enumerate}
\hfill $\blacksquare$

We now come to the case of slice regular functions.
In the sequel, we will denote by $W$ the operator  acting on function defined from  $\Lambda$ to $\mathbb C$ as follows:
$$W[f](z)= \overline{f(\bar z)},\quad \forall z\in\Lambda,$$
for any complex valued function $f$ with domain  $\Lambda$.

Let now $\Omega\subset\mathbb H$ be an axially symmetric s-domain. For any ${\bf i}\in\mathbb S^2$ denote $\Omega_{\bf i}^+ :=\{  x+{\bf i} y \in \Omega_{\bf i}
\mid \   y > 0   \}$, {}
$\Omega_{\bf i}^- :=\{ x+{\bf i } y \in \Omega_{\bf i}
\mid \   y < 0   \}$, and  {} $\Omega_0 := \Omega\cap \mathbb R$.

\begin{defin}  \label{def21} Let $\Omega$  be an axially symmetric s-domain.
  \begin{enumerate}\item  We define $\displaystyle \widetilde{Hol} (\Omega^+_{\bf i})= \{  f\in Hol(\Omega^+_{\bf i}) \ \mid \   \textrm{  for any  } x\in \Omega_0  \  \textrm{there exists}  \lim_{\Omega^+_{\bf i}\ni z \to x } f(z) \in\mathbb R   \} .$
\item  For any $f\in \widetilde{Hol} (\Omega^+_{\bf i})$ we define the function on $\Omega_{\bf i}\subseteq\mathbb C_{\bf i}$ by  $$E_{\bf i }[f](z)=\left\{ \begin{array}{ll}  f(z) &  \textrm{ if  }  z\in \Omega^+_{\bf i}\\ \\ \displaystyle  \lim_{\Omega^+\ni w \to z } f(w) \in\mathbb R   & \textrm{ if } z\in \Omega_0  \\ \\ \overline{f(\bar z)} & \textrm{ if } z\in \Omega^-_{\bf i}  \end{array}  \right.$$

\item Finally, let $$E_{\bf i}( \ \widetilde{Hol}(\Omega^+_{\bf i})\  )=\{ E_{\bf i}[f] \ \mid \ f\in \widetilde{Hol} (\Omega^+_{\bf i}) \}.$$\end{enumerate}
\end{defin}

According to the notation of the  previous definition,  the Schwarz reflection principle implies  that
$Hol_c(\Omega_{\bf i}) = E_{\bf i}( \ \widetilde{Hol}(\Omega^+_{\bf i})\  )$, and as  $Hol(\Omega_{\bf i})= Hol_c(\Omega_{\bf i})  + Hol_c(\Omega_{\bf i} ){\bf i}$, then
$$Hol(\Omega_{\bf i})= E_{\bf i}( \ \widetilde{Hol}(\Omega^+_{\bf i})\  ) + E_{\bf i}( \ \widetilde{Hol}(\Omega^+_{\bf i})\  ){\bf i}.$$
Note that
$Hol_c(\Omega_{\bf i})$ and $ \widetilde{Hol}(\Omega^+_{\bf i})$ are linear spaces over $\mathbb R$.

\begin{ob}\label{rem1} {\rm Let ${\bf i},{\bf j}\in\mathbb S^2$ be such that  ${ {\bf j} }  \perp {\bf i}$.
    Then
    $$\mathcal{SR}(\Omega)=  P_{\bf i}\circ  E_{\bf i}[ \ \widetilde{Hol}(\Omega^+_{\bf i})\  ] +  P_{\bf i}\circ E_{\bf i}[ \ \widetilde{Hol}(\Omega^+_{\bf i})\  ]{\bf i}   + P_{\bf i}\circ  E_{\bf i}[ \ \widetilde{Hol}(\Omega^+_{\bf i})\  ] {\bf j}+  P_{\bf i}\circ E_{\bf i}[ \ \widetilde{Hol}(\Omega^+_{\bf i})\  ]{\bf i}{\bf j} .  $$
Thus any slice regular function can be written in terms of four elements, each of them belonging $\widetilde{Hol}(\Omega^+_{\bf i})$.
By changing $\Omega_{\bf i}^+$ with $\Omega_{\bf i}^-$  we obtain a decomposition of $\mathcal{SR}(\Omega)$  in terms of $\widetilde{Hol}(\Omega^-_{\bf i})$.
}
\end{ob}

\begin{pro}\label{pro2}
Let
$ f\in \widetilde{Hol}(\Omega^+_{\bf i})$,  then its extension $P_{\bf i}[f]$ to the whole $\Omega$ is given by
$$ Re f (q_0 + {\bf i} |\underline{q}| )  \ +  \ I_q \  Im  f (q_0 + {\bf i} |\underline{q}| ), $$
for all $q\in \Omega$.
\end{pro}

Proof. Formula (\ref{extPi}) immediately gives
  $$
   P_{\bf i}[f](q)=P_{\bf i}[E_{\bf i}[f]](q_0 +  I_q |\underline{q}| )=\frac{1}{2}\left[(1+ I_q{\bf i})\overline{f(q_0 +|\underline{q}|{\bf i})} + (1- I_q {\bf i}) f(q_0 +|\underline{q}| {\bf i})\right]
   $$
$$= {\rm Re} f (q_0 + {\bf i} |\underline{q}| )  \ +   I_q \,  {\rm Im}  f (q_0 + {\bf i} |\underline{q}| ).$$
\hfill $\blacksquare$

\begin{pro}
Let
$ f\in \mathcal {SR}(\Omega)$,  then there exist  conjugated harmonic functions $u_n,v_n \in C^2(\Omega_{\bf i}^+,\mathbb R)$,  $n=0,1,2,3$, with
 $\displaystyle \lim_{\Omega^+\ni w \to z }v_{n }(w)=0  $,  for each $ z\in\Omega_0    $,    such that
$$f (q) =f(q_0 + I_q |\underline{q}| )= \sum_{n=0}^3 u_n (q_0 + {\bf i} |\underline{q}| )e_n  \ +  \ I_q \ \sum_{n=0}^3 v_n (q_0 + {\bf i} |\underline{q}| )e_n, $$
for all $q\in \Omega$, where $e_0=1, \  e_1={\bf i}, \ e_2={\bf j}, \ e_3={\bf ij}$.
\end{pro}
Proof.
It is a direct consequence of Remark \ref{rem1} and Proposition \ref{pro2}.
\hfill $\blacksquare$
\vskip 0.3truecm
From the previous proposition one has that:
\begin{enumerate}
\item If $  q=q_0 \in\Omega_0$ then $\displaystyle f (q_0) = \sum_{n=0}^3 u_n (q_0 )e_n .$
\item $\displaystyle  f (\bar q) = \sum_{n=0}^3 u_n (q_0 + {\bf i} |\underline{q}| )e_n  \, -  I_q \, \sum_{n=0}^3 v_n (q_0 + {\bf i} |\underline{q}| )e_n.$
\end{enumerate}

\begin{defin}Let ${\bf i} \in \mathbb S^2$. We set:
\begin{enumerate}
\item $\widetilde{\mathcal A}(\Omega_{\bf i}^+) =\widetilde{Hol} (\Omega_{\bf i}^+) \cap \mathcal L_2(\Omega_{\bf i}^+) $
\item $\mathcal A(\Omega_{\bf i}) =\mathcal{SR} (\Omega) \cap \mathcal L_2(\Omega_{\bf i}) $.
\end{enumerate}
\end{defin}
The following result shows that the Bergman type space $\mathcal A(\Omega_{\bf i})$ can be written in terms of the Bergman type spaces on half of a slice $\Omega_{\bf i}^+$.
Moreover, the inner product of two elements can be computed on half of a slice.
\begin{pro} The following facts hold:\begin{enumerate}
\item
\begin{equation}\label{decomp}\mathcal A(\Omega_{\bf i}) = P_{\bf i}\circ  E( \ \widetilde{\mathcal A}(\Omega_{\bf i}^+)  \  ) +P_{\bf i}\circ  E( \ \widetilde{\mathcal A}(\Omega_{\bf i}^+) \  ){\bf i}+P_{\bf i}\circ E( \ \widetilde{\mathcal A}(\Omega_{\bf i}^+)  \  ) {\bf j}+P_{\bf i}\circ  E( \ \widetilde{\mathcal A}(\Omega_{\bf i}^+) \  ){\bf i}{\bf j}.
\end{equation}

  \item Let $f,g \in \mathcal A(\Omega_{\bf i})$   and let  $f_1f_2,g_1,g_2 \in  Hol (\Omega_{\bf i})\cap \mathcal L_2(\Omega_{\bf i},\mathbb C({\bf i}) )$  be  such that  $f=P_{\bf i}[f_1] + P_{\bf i}[f_2]{\bf j}$  and $g=P_{\bf i}[g_1] + P_{\bf i}[g_2]{\bf j}  $.  Then
  $$\displaystyle  <f, g>_{\mathcal A(\Omega_{\bf i})} =  \left(  \  \displaystyle  <f_1\mid_{\Omega_{\bf i}^+}, g_1\mid_{\Omega_{\bf i}^+}>_{\mathcal L(\Omega^+_{\bf i},\mathbb C({\bf i}))}
+ <W[f_1]\mid_{\Omega_{\bf i}^+}, W[g_1]\mid_{\Omega_{\bf i}^+}>_{\mathcal L(\Omega^+_{\bf i},\mathbb C({\bf i}))}
+\right.$$
$$  \displaystyle +  \left. \overline{ <f_2\mid_{\Omega_{\bf i}^+}, g_2\mid_{\Omega_{\bf i}^+}>}_{\mathcal L(\Omega^+_{\bf i},\mathbb C({\bf i}))}
+    \overline{ <W[f_2]\mid_{\Omega_{\bf i}^+}, W[g_2]\mid_{\Omega_{\bf i}^+}>}_{\mathcal L(\Omega^+_{\bf i},\mathbb C({\bf i}))}  \    \right) +$$
$$ + \displaystyle  \left( \
<f_1\mid_{\Omega_{\bf i}^+}, g_2\mid_{\Omega_{\bf i}^+}>_{\mathcal L(\Omega^+_{\bf i},\mathbb C({\bf i}))}
+ <W[f_1]\mid_{\Omega_{\bf i}^+}, W[g_2]\mid_{\Omega_{\bf i}^+}>_{\mathcal L(\Omega^+_{\bf i},\mathbb C({\bf i}))}
 \right. -  $$
$$\displaystyle \left.  - \overline{   <f_2\mid_{\Omega_{\bf i}^+}, g_1\mid_{\Omega_{\bf i}^+}>}_{\mathcal L(\Omega^+_{\bf i},\mathbb C({\bf i}))}
+  \overline{ <W[f_2]\mid_{\Omega_{\bf i}^+}, W[g_1]\mid_{\Omega_{\bf i}^+}>}_{\mathcal L(\Omega^+_{\bf i},\mathbb C({\bf i}))}
    \   \right) {\bf j}$$

  \item Let $f\in \mathcal A(\Omega_{\bf i})$. Then
  $$\|f\|_{\mathcal A(\Omega_{\bf i})}^2 =  2  \left[  \int_{\Omega^+} |f\mid_{\Omega_{\bf i}}  ^+|^2 d\sigma\right]
 $$
 \end{enumerate} \end{pro}
 Proof.
\begin{enumerate}
\item The decomposition \eqref{decomp} is a  consequence of Remark \ref{rem1}.
\item Let us compute the inner product $<f, g>_{\mathcal A(\Omega_{\bf i})}$. We have:
$$\displaystyle  <f, g>_{\mathcal A(\Omega_{\bf i})} =   \int_{\Omega_{\bf i}}  \overline{( f_1+ f_2{\bf j}  ) } (g_1+g_2{\bf j})   d\sigma=$$
  $$\displaystyle   =   \int_{\Omega_{\bf i}}  \overline{ f_1 }g_1  d\sigma - {\bf j}
 \int_{\Omega_{\bf i}}  \overline{f_2 }g_2   d\sigma{\bf j}
+ \int_{\Omega_{\bf i}}  \overline{ f_1 }g_2   d\sigma{\bf j} - {\bf j   } \int_{\Omega_{\bf i}}  \overline{f_2 } g_1  d\sigma=$$

$$=\displaystyle  <f_1, g_1>_{\mathcal L(\Omega_{\bf i},\mathbb C({\bf i}))}  -{ \bf j}
<f_2, g_2>_{\mathcal L(\Omega_{\bf i},\mathbb C({\bf i}))}  {\bf j} +
<f_1, g_2>_{\mathcal L(\Omega_{\bf i},\mathbb C({\bf i}))}  {\bf j} -
{\bf j } <f_2, g_1>_{\mathcal L(\Omega_{\bf i},\mathbb C({\bf i}))}=$$

$$=\left(  \  \displaystyle  <f_1, g_1>_{\mathcal L(\Omega_{\bf i},\mathbb C({\bf i}))} + \overline{ <f_2, g_2>}_{\mathcal L(\Omega_{\bf i},\mathbb C({\bf i}))}   \    \right) +  \left( \
   <f_1, g_2>_{\mathcal L(\Omega_{\bf i},\mathbb C({\bf i}))}  - \overline{   <f_2, g_1>}_{\mathcal L(\Omega_{\bf i},\mathbb C({\bf i}))}     \   \right) {\bf j}$$

$$=\left(  \  \displaystyle  <f_1\mid_{\Omega_{\bf i}^+}, g_1\mid_{\Omega_{\bf i}^+}>_{\mathcal L(\Omega^+_{\bf i},\mathbb C({\bf i}))}
+ <W[f_1]\mid_{\Omega_{\bf i}^+}, W[g_1]\mid_{\Omega_{\bf i}^+}>_{\mathcal L(\Omega^+_{\bf i},\mathbb C({\bf i}))}
+\right.$$
$$  \displaystyle +  \left. \overline{ <f_2\mid_{\Omega_{\bf i}^+}, g_2\mid_{\Omega_{\bf i}^+}>}_{\mathcal L(\Omega^+_{\bf i},\mathbb C({\bf i}))}
+    \overline{ <W[f_2]\mid_{\Omega_{\bf i}^+}, W[g_2]\mid_{\Omega_{\bf i}^+}>}_{\mathcal L(\Omega^+_{\bf i},\mathbb C({\bf i}))}  \    \right) +$$
$$ + \displaystyle  \left( \
<f_1\mid_{\Omega_{\bf i}^+}, g_2\mid_{\Omega_{\bf i}^+}>_{\mathcal L(\Omega^+_{\bf i},\mathbb C({\bf i}))}
+ <W[f_1]\mid_{\Omega_{\bf i}^+}, W[g_2]\mid_{\Omega_{\bf i}^+}>_{\mathcal L(\Omega^+_{\bf i},\mathbb C({\bf i}))}
 \right. -  $$
$$\displaystyle \left.  - \overline{   <f_2\mid_{\Omega_{\bf i}^+}, g_1\mid_{\Omega_{\bf i}^+}>}_{\mathcal L(\Omega^+_{\bf i},\mathbb C({\bf i}))}
 - \overline{ <W[f_2]\mid_{\Omega_{\bf i}^+}, W[g_1]\mid_{\Omega_{\bf i}^+}>}_{\mathcal L(\Omega^+_{\bf i},\mathbb C({\bf i}))}
   \   \right) {\bf j}$$
\item
$$\|f\|_{\mathcal A(\Omega_{\bf i})}^2 =
\|f_1\|_{\mathcal A(\Omega_{\bf i})}^2   +\|f_2\|_{\mathcal A(\Omega_{\bf i})}^2
$$
$$=
  2\left[     \int_{\Omega^+} |f_1\mid_{\Omega_{\bf i}}  ^+|^2 d\sigma
+    \int_{\Omega^+_2} | f\mid_{\Omega_{\bf i}}  ^+|^2 d\sigma \right]
$$
$$
=  2  \left[  \int_{\Omega^+} |f\mid_{\Omega_{\bf i}}  ^+|^2 d\sigma\right]$$

\end{enumerate}
\hfill $\blacksquare$

\end{document}